\documentclass{amsart}
\usepackage{amscd,amssymb}

\theoremstyle{plain}
\newtheorem{thm}{Theorem}
\newtheorem{lem}[thm]{Lemma}
\newtheorem{definicja}[thm]{Definition}

\theoremstyle{definition}

\theoremstyle{remark}

\numberwithin{equation}{section}

\begin{document} 
\textheight = 650pt 

\title{ Polish group actions and effectivity} 
\author{ Barbara Majcher-Iwanow  } 

\maketitle 

\begin{quote}
{\bf Abstract.} We extend the result of Nadel describing the relationship
between approximations of canonical Scott sentences and admissible sets
to the general case of orbit equivalence relations induced on an arbitrary
Polish space by a Polish group action.

\bigskip


{\em Classification:} 03E15, 03C70
\end{quote}

\bigskip

\section{Introduction}

In the paper we extend the following result of Nadel (see \cite{nadel})
to the general case of Polish $G$-spaces.
\begin{quote}
{\bf Theorem}(Nadel) {\em Let $L$ be a countable language,
${\mathbb{A}}$ be an admissible set
and ${\frak{M}}$ be an $L$-structure in ${\mathbb{A}}$.
Then for any $L$-structure ${\frak{N}}$, if ${\frak{M}}$ and
${\frak{N}}$ satisfy the same sentences from the admissible fragment
$L_{\mathbb{A}}$, then they satisfy the same sentences
of quantifier rank $\alpha\le {\it{o}}({\mathbb{A}}).$ }
\end{quote}
The theorem can be formulated in terms of the logic action as follows.
Consider a countable relational language $L= (R^{n_i}_i)_{i\in I}$.
There is an obvious one-to-one correspondence
${\mathfrak{M}}\to x_{\mathfrak{M}}$ between the set of
all countable $L$-structures and  the set
$X_L=\prod_{i\in I}2^{\omega^{n_i}}$.
The set $X_L$ equipped with the product topology becomes a Polish space,
the space of all $L$-structures on $\omega$
(see Section 2.5 in \cite{bk} or Section 2.D of \cite{becker}
for details).
The group $S_{\infty}$ of all permutations of $\omega$ has
the natural continuous action on the space $X_L$.
It is called the {\em logic action} of $S_{\infty}$ on $X_L$.
Given a structure $x\in X_L$, the orbit of $x$ under the logic action
consists of all structures isomorphic to $x$.
Now the Nadel's theorem reads as follows.
\begin{quote}{\em Let ${\mathbb{A}}$ be an admissible set
and $L\in {\mathbb{A}}$ be a countable relational language.
Let $x,y\in X_L$ and $x\in {\mathbb{A}}$.
Then if $x$ and $y$ are in the same invariant Borel
sets of the form $B_{\sigma}=\{z\in X_L: z\models \sigma\}$
for some $\sigma\in L_{\mathbb{A}}$,
then for every $\alpha\le o({\mathbb{A}})$,
$x$ and $y$ are in the same invariant $\Pi_{\alpha}^0$-subsets of $X_L$.}
\end{quote}
The theorem can be viewed as an assertion to what extent a structure
$x$ of a given language $L$ is determined (up to isomorphism)
by the set of sentences from $L_{\mathbb{A}}$ true in $x$.
We shall ask the general question:
\begin{quote}{\em Given a Polish $G$-space $X$ which
is in some sense coded in an admissible set ${\mathbb{A}}$ and $x\in X$,
to what extent is the orbit $G\cdot x$ determined
by the family of those invariant Borel sets
containing $x$ which have Borel codes in ${\mathbb{A}}$}?
\end{quote}
In the paper we give an answer to this question.
It generalizes the Nadel's theorem for any Polish group $G$
and any Polish $G$-space $X$.
It is formulated in the second part of Section 2.
In the first part of this section we  settle the parallel question
on the ground of effective descriptive set theory.
This is the main result of the paper.
A simple proof of this theorem applies some tools
and facts from the generalized
Scott analysis of continuous actions of Polish groups
on Polish spaces developed by Hjorth in \cite{hjorth}.

\section{Preliminaries}

In the first part of this section we recall standard notation and
facts concerning Polish group actions as well as the brief description
of Hjorth method.
In the second one we give a summary introduction
to effective descriptive set theory.

\subsection{Notation}
\parindent1em
A {\em Polish space} ({\em group}) is a separable, completely
metrizable topological space (group).
We shall write $(X,d)$ if $d$ is a compatible metric for $X$.
If a Polish group $G$ continuously acts on a Polish space $X$,
then we say that $X$ is a {\em Polish $G$-space}.
We  say that a subset of $X$ is {\em invariant} if
it is $G$-invariant.
All basic facts concerning Polish $G$-spaces can
be found in \cite{bk}, \cite{hjorth} and \cite{kechris}.
\parskip0pt

The Vaught $*$-transform of a set $B\subseteq X$
with respect to an open $H\subseteq G$ is the set
$B^{*H}=\{ x\in X:\{ g\in H:gx\in B\}\mbox{ is comeagre in } H\}$.
A set $B$ is invariant if and only if $B^{*G}=B$.
\bigskip

\parindent0em
{\bf Definition A} (Hjorth) {\em Let $x,y\in X$.
For basic open $V, W\subseteq G$ we define the relation
$(y,V)\le_{\alpha} (x,W)$ by simultaneous induction on the ordinal $\alpha$.}

$$\begin{array}{l}(y,V)\le_1 (x,W)\  \mbox{ {\em if} }
\ \overline{V\cdot y}\subseteq \overline{W\cdot x}\\
\\
(y,V)\le_{\alpha+1}(x,W) \mbox{ {\em if for every basic open} } V'\subseteq V\\
\mbox{ {\em there is a basic open} } W'\subseteq W
\mbox{ {\em such that} } (x,W')\le_{\alpha}(y,V')
\\
\\
(y,V)\le_{\lambda}(x,W)\mbox{ {\em if} } (y,V)\le_{\alpha} (x,W)
\mbox{ {\em for every} }\alpha<\lambda \\
\mbox{ {\em if} }\lambda \mbox{ {\em is limit} }.\end{array}$$

\parindent1em
It is shown in \cite{hjorth} that the relation $\le_{\alpha}$ is transitive
and $(y,V)\le_{\alpha}(x,W)$ implies $(y,V)\le_{\beta}(x,W)$,
whenever $\beta\le\alpha$.
Moreover for every ordinal $\alpha$ and $x\in X$ we have
$(x,G)\le_{\alpha}(x',G)$ whenever $x'\in G\cdot x$.
For our purpose the most important property is as follows.
\bigskip

\parindent0em
{\bf Lemma B} (Hjorth)
{\em
Let $x,y\in X$, $\ V,W\subseteq G$ be basic open sets
and $\alpha$ be a countable ordinal such that
$(y,V)\le_{\alpha}(x,W)$.
Then for every ${\bf \Pi}_{\alpha}^0$-set $B\subseteq X$
if $x\in B^{*W}$, then $y\in B^{*V}$.} 

\subsection{Review of basic notions and facts from
effective descriptive set theory.}

\parindent1em
We assume here acquaintance with rudiments of recursion theory.
We recall only those notions and facts we shall use in the paper.
The review below is based on  \cite{MaKe} and \cite{mos} (Chapters 3,7),
where one can find a complete exposition.

We shall consider merely {\em recursively presented} Polish spaces.
A Polish space $(X,d)$ is recursively presented if it is endowed
with a {\em recursive presentation},
i.e. a sequence $\{r_i:i\in\omega\}$
densely contained in $X$ and such that the $(i,j,k,m)$-relations
$\ d(r_i,r_j)\le m/k+1$ and  $d(r_i,r_j)<m/k+1\ $ are recursive.
The class of {\em recursively presented} spaces includes
$\omega$, the Baire space ${\mathcal{N}}$, the Cantor space ${\mathcal{C}}$,
the reals ${\mathcal{R}}$ and it is closed under finite products.
\parskip0pt

If $X$ is a recursively presented space
with a fixed recursive presentation,
then one can naturally define the canonical basis
of open neighbourhoods of $X$, an effective enumeration
$\{U_n^X:n\in\omega\}$ of the basis
and a ternary recursive function $g(k,l,m)$ with $Dom(g)=\omega^3$ so that
$$
U^X_k \cap U^X_l = \bigcup_{m\in\omega} U^X_{g(k,l,m)}.
$$
We call a set $A\subseteq X$ {\em semirecursive}
or {\em effectively open} if there is a recursive
function $\ s:\omega\to\omega\ $ such that  $\ A=\bigcup\limits_nU^X_{s(n)}$.
\parskip0pt

The pointclass of semirecursive pointsets
(of recursively presented spaces) is denoted by
$\Sigma^0_1$ and $\Pi^0_1$ stand for the pointclass
of {\em effectively closed} (i.e. complements of semirecursive) pointsets.
The ambiguous class $\Delta^0_1=\Sigma^0_1\cap \Pi^0_1$
\footnote{notice using lightface font for distinguishing
effective classes from their non-effective analogues
printed later on in boldface}
consists of {\em recursive} pointsets.
The pointclass $\Delta^0_1$ contains the empty set,
every product space, every recursive relation on $\omega^k$,
every basic neighbourhood $U_n^X$ and the basic neighbourhood
relation $\ \{(x,n): x\in U_n^X\}$ for each $X$.\parskip0pt

If we start from the pointclass $\Sigma^0_1$ of effectively open sets
and $\omega^{CK}_1$-times alternately take the operation of complementation
and the operation of "effective countable unions"
(the Church-Kleene ordinal $\omega_1^{CK}$
is the first non-recursive ordinal), then we obtain the pointclass
$$
HYP=\bigcup\limits_{\zeta<\omega_1^{CK}}\Sigma_{\zeta}^0
$$
of {\em hyperarithmetic} pointsets,
the effective analogue of the pointclass of Borel pointsets.
On the other hand if we apply to $\Sigma^0_1$ the operations
of complementation and projection along Baire space ${\mathcal{N}}$,
then we build the pointclasses
$$\Sigma^1_1=\exists^{\mathcal{N}}\Pi^0_1,\
\Pi^1_1=\neg\Sigma^1_1 \ \mbox{ and }
\ \Delta^1_1=\Sigma^1_1\cap\Pi^1_1.$$
The following statement is true for every recursively presented space $X$.
\bigskip

\parindent0em
{\bf Theorem C} (The Suslin-Kleene theorem)
{\em For every $A\subseteq X$ we have}
$$A\in HYP\ \mbox{ iff }\ A\in \Delta^1_1.$$

\parindent0em
{\bf Definition D}
{\em Let $X,Y$ be  recursively presented spaces with bases
$\{U^X_n\}$ and $\{U^Y_n\}$ respectively,
$\Gamma=\Sigma^0_{\xi}, \ \Sigma^1_1\mbox{ or } \Delta^1_1$.}\parskip2pt

(1) {\em A function $f:X\to Y$ is
$\Gamma$-recursive if } $\{(x,n): f(x)\in U^Y_n\}\ \in \ \Gamma.$

(2) {\em A partial function $f:X\to Y$ is
$\Gamma$-recursive on its domain if there is some

$P\subseteq X\times \omega$
such that $P\in\Gamma$ and }
$\{(x,n): f(x)\in U^Y_n\}\ =\ P\cap (dom(f)\times\omega).$

(3) {\em An element $x\in X$ is
$\Gamma$-recursive if }
$\{n:x\in U^X_n\}\in \ \Gamma.$
\bigskip

\parindent1em
$\Sigma^0_1$-recursive functions (points) are simply called {\em recursive}
and $\Delta^1_1$-recursive - {\em hyperarithmetic}.
The pointclass $\Delta_1^1$ satisfies the following closure properties.
\bigskip

\parindent0em
{\bf Theorem E} {\em The pointclass $\Delta_1^1$ is closed under
$\vee,\wedge$, negation, existential and universal quantification
over $\omega$ and substitution of $\Delta^1_1$-recursive functions}
\footnote{This means that for every partial $\Delta^1_1$-recursive function
$f:X\to Y$ and every $\Delta^1_1$-set $B\subseteq Y$ there is
a $\Delta^1_1$-set $A\subseteq X$ such that $f^{-1}[B]=A\cap dom(f)$.}.
\bigskip

\parindent1em
For every space $Z$, an element $z\in Z$ and a pointclass
$\Gamma=\Sigma^0_{\xi}, \ \Pi^0_{\xi},\ \Sigma^1_1\mbox{ or } \Pi^1_1$
we shall consider also the relativized  pointclass
$\Gamma(z)$: a set $A\subseteq X$ is in
$\Gamma(z)$ if there is some $Q\subseteq Z\times X$
such that $Q\in \Gamma$  and
$A=Q_z$, where $Q_z=\{x\in X: (z,x)\in Q\}$.
Then we define the ambiguous classes
$\Delta^0_1(z)=\Sigma^0_1(z)\cap \Pi^0_1(z)$ and
$\Delta^1_1(z)=\Sigma^1_1(z)\cap \Pi^1_1(z)$.
The sets in $\Sigma^0_1(z)$ are called {\em semirecursive in} $z$,
the sets in $\Delta^1_1(z)$ - {\em hyperarithmetic in} $z$.

In the obvious way we may relativize the notions of recursive and
hyperarithmetic functions.
It should be mentioned  that the relativized versions of the theorems cited
above remain true.\parskip0pt

It is easy to see that a subset of a given space $X$ is open
if and only if it is semirecursive in some $z\in {\mathcal{N}}$.
In fact for every recursively presented space $X$ there is a semirecursive
$P\subseteq  {\mathcal{N}}\times X$
(called a {\em good parametrization system})  which is universal
for the pointclass of open subsets of $X$ so that
for every $A\subseteq X$ we have
$$A\in \Sigma^0_1 \ \mbox{ iff } \ A=P_{\epsilon}
\ \mbox{ for a recursive } \epsilon\in {\mathcal{N}}.
\footnote{This is a part of The Good Parametrization Lemma,
see \cite{mos} 3H.1}$$

Starting with such a good parametrization system
$P\subseteq {\mathcal{N}}\times\omega\times\omega$ for
subsets of $\omega\times\omega$ we can define some special
coding of Borel sets by elements of ${\mathcal{N}}$.
This method is described in detail in \cite{mos} (Chapter 7.A),
we shall only recall some notation and properties

Using $P$ we define by recursion on the countable ordinal $\xi$
an increasing family $\{BC_{\xi}\}_{\xi}$ of subset of
${\mathcal{N}}$ with the union
$BC=\bigcup_{\xi}BC_{\xi}$ - the set of Borel codes.

Then for every recursively presented space $X$ we can define
a coding function $\pi:BC\to Borel(X)$ so that
$\pi|_{BC_{\xi}}:BC_{\xi}\to {\bf\Sigma}^0_{\xi}(X)$ is "onto",
for every countable $\xi$.
Given $A\in Borel(X)$ we say that an irrational $\alpha\in BC$
is a {\em Borel code} of $A$ whenever $\pi(\alpha)=A$.
This coding function has the following important property.
\bigskip

\parindent0em{\bf Lemma F}
{\em Let $z\in {\mathcal{N}}$ and $A\subseteq X$.
Then $A$ is hyperarithmetic in $z$ if and only if it has
a Borel code recursive in $z$}.\parskip0pt

\section{General versions of the Nadel's theorem}

\parindent1em
The section is divided into two parts.
In the first one we shall prove an effective counterpart 
of the Nadel's theorem.
In the second part we translate
these results into the language of admissible sets.

\subsection{Hyperarithmetic actions and Borel sets}

From now on we shall always assume that a Polish group
$G$ is a recursively presented space
and  $\{\hat{g}_i:i\in\omega\}$
is a fixed recursive presentation of $G$.
Let $\{V_i:i\in\omega\}$ be an effective enumeration
of the canonical basis of $G$.
We additionally assume that $V_0=G$.
A Polish $G$-space $X$ is recursively presented and $\{U_i:i\in\omega\}$
is an effective enumeration of the canonical basis of $X$.

It is proved in \cite{hjorth} that the relation
$\{(x,y,n,k): (x,V_k)\le_{\xi}(y,V_n)\}\subseteq X^2\times\omega^2$
is Borel for every countable ordinal $\xi$.
The following lemma is an effective version of this statement.

\begin{lem}\label{hypH}
Let $G$, $X$ be recursively presented, the group operations
of $G$ and the $G$-action on $X$ be hyperarithmetic.
Then for every $\xi<\omega_1^{CK}$ the set
$$H_{\xi}=\{(x,y,n,k): (x,V_k)\le_{\xi}(y,V_n)\}$$
is hyperarithmetic.
\end{lem}

{\bf Proof.}
First consider $H_1=\{(y,x,k,n): \overline{V_ky}\subseteq \overline{V_nx}\}$.
We have
$$\begin{array}{c}
(y,x,k,n)\in H_1\\
\Updownarrow\\
(\forall i)(y\in V^{-1}_kU_i\Rightarrow x\in V^{-1}_nU_i)\\
\Updownarrow \\
(\forall i)\Big((y,x)\in [(X\setminus V^{-1}_kU_i)\times X]\cup
[V^{-1}_kU_i\times V^{-1}_nU_i]\Big)\\
\Updownarrow\\
(y,x,k,n)\in \bigcap_i\Big([(X\setminus V^{-1}_kU_i)\times X]\cup
[V^{-1}_kU_i\times V^{-1}_nU_i]\Big)\times\{k\}\times\{n\}.
\end{array}$$
Hence $H_1=\bigcup_{k,n}\bigcap_i\Big([(X\setminus V^{-1}_kU_i)\times X]
\cup [V^{-1}_kU_i\times V^{-1}_nU_i]\Big)\times\{k\}\times\{n\}$.\bigskip

Now observe that if $V\subseteq G$ and $U\subseteq X$ are open,
then $VU=\bigcup\{\hat{g}_jU: \hat{g}_j\in V, j\in\omega\}$.
Indeed, take any $g\in V$ and $x\in U$.
It follows from the continuity of the group operations of $G$ and the
$G$-action that there are open
$g\in W\subseteq V$ and $x\in A\subseteq U$ such that
$WW^{-1}W\subseteq V$ and $W^{-1}WA\subseteq U$.
Then for any $f\in W$ we  have $g\in fW^{-1}W$ and
so $gx\in fW^{-1}WA\subseteq fU$.
In particular $V_k^{-1}U_i=\bigcup\{\hat{g}_j^{-1}U_i:\hat{g}_j\in V_k\}$.
Since the $G$-action is hyperarithmetic,
every element $\hat{g}_j$ of the recursive presentation
$\{\hat{g}_i:i\in\omega\}$ is recursive and every basic open $U_i$
is recursive, then by closeness of $\Delta^1_1$ pointclass under
substitution property of $\Delta^1_1$-functions, $\hat{g}_j^{-1}U_i\in HYP$.
Hence by other closure properties from Theorem E,
the set $H_1$ is hyperarithmetic.\parskip0pt

Next we have
$$\begin{array}{c}
(y,x,k,n)\in H_{\mu+1}\\
\Updownarrow\\
(\forall r)(\exists s)\Big(V_r\subseteq V_k\Rightarrow
(V_s\subseteq V_n)\wedge (x,V_s)\le_{\mu}(y,V_r)\Big)\\
\Updownarrow \\
(\forall r)(\exists s)\Big(V_r\subseteq V_k\Rightarrow
(V_s\subseteq V_n)\wedge (x,y,s,r)\in H_{\mu})\Big)\\
\Updownarrow\\
(x,y)\in \bigcap\limits_{r\atop{V_r\subseteq V_k}}\Big(
\bigcup\limits_{s\atop{V_s\subseteq V_n}}\pi_{1,2}
\Big(H_{\mu}\cap (X^2\times\{r\}\times\{s\})\Big)\Big)\\
\Updownarrow \\
(y,x,k,n)\in s\Big(\bigcap\limits_{r\atop{V_r\subseteq V_k}}\Big(
\bigcup\limits_{s\atop{V_s\subseteq V_n}}\pi_{1,2}
\Big(H_{\mu}\cap (X^2\times\{r\}\times\{s\})\Big)\Big)\Big)
\times\{k\}\times\{n\},
\end{array}$$
where $\pi_{1,2}(x,y,r,s)=(x,y)$ and $s(x,y)=(y,x)$.

Since $G$ is recursively presented, the relation
$\{(k,r): V_r\subseteq V_k\}$ is recursive.
Hence using closure properties of $\Delta^1_1$ we argue that
if $H_{\mu}$ is hyperarithmetic, then so is $H_{\mu+1}$.\parskip0pt

On the other hand if $\xi$ is a recursive limit ordinal and $H_{\mu}$
is hyperarithmetic for every $\mu<\xi$,
then $H_{\xi}=\bigcap_{\mu<\xi}H_{\mu}$
is an "effective intersection" of hyperarithmetic sets -
thus is also hyperarithmetic.
$\Box$\bigskip

Now we are ready to prove the main result of this part.
\begin{thm}\label{nadef}
Assume that $G$, $X$ are  presented recursively in $\zeta$,
the group operations of $G$ and the $G$-action on $X$
are hyperarithmetic in $\zeta$.
If $x,y\in X$, $x$ is hyperarithmetic in $\zeta$ and $x,y$
are in the same invariant sets hyperarithmetic in $\zeta$,
then for every $\alpha\le\omega_1^{CK}(\zeta)\ $
$x,y$ are in the same invariant Borel sets of Borel rank $\alpha$.

\end{thm}

{\bf Proof.} Assume that $\zeta\in \Delta^0_1$.
Let $\alpha<\omega^{CK}_1$ and  $B\subseteq X$
be an invariant ${\bf{\Pi}}^0_{\alpha}$-set containing $x$.
By the lemma above $H_{\alpha}$ is hyperarithmetic.
Since $x$ is hyperarythmetic, then the set
$\{z\in X:(z,G)\le_{\alpha}(x,G)\}=
\pi_1(H_{\alpha}\cap (X\times \{x\}\times\{0\}\times\{0\})$
is an invariant hyperarithmetic set containing $x$.
Thus $y$ is an element of this set, i.e.
$(y,G)\le_{\alpha}(x,G)$.
This by Lemma B implies $y\in B$.\parskip0pt

Now consider the case $B\in {\bf \Pi}^0_{\omega^{CK}_1}$.
Since ${\omega^{CK}_1}$ is a limit ordinal,
then $B$ is an intersection of ${\omega^{CK}_1}$ invariant sets
of Borel rank $\xi< {\omega^{CK}_1}$.
If $x\in B$, then $x$ belongs to each of the sets.
Hence by the first part of the proof $y\in B$.\parskip0pt

In the case of arbitrary $\zeta$ we use the same arguments
based on the obvious relativization of Lemma \ref{hypH}.
$\Box$

\subsection{General version of the Nadel's theorem}
We assume that the reader is familiar with
the most basic notions of \emph{admissible sets}. Any necessary background
can be easily provided by \cite{barwise} and \cite{ershov}. \parskip0pt

We only remind the reader that an admissible set ${\mathbb{A}}$ is a
transitive model of KPU, in the sense of \cite{barwise}. Such models are
considered as two-sorted structures of some language $L$ with symbols
$\emptyset ,\in$, where one of the sorts corresponds
to \emph{urelements} and
usually forms a relational first-order structure with respect to the symbols
of $L$ distinct from $\emptyset$ and $\in$. Here we assume that
${\mathbb{A}}$ satisfies KPU with respect to all formulas of $L$
(${\mathbb{A}}$ is admissible with respect to $L$ in terms of \cite{nadel}).
\parskip0pt

As we recalled in Section 1.2 for every recursively presented space $X$
we can define a partial function
$\pi: {\mathcal{N}}\ {\stackrel{onto}\rightarrow}\ Borel(X)$ with
domain the set $BC$ of Borel codes,
so that a set $B\subseteq X$ is hyperarithmetic
(resp. hyperarithmetic in $z$) if and only if it has a recursive
(resp. recursive in $z$) Borel code.
Thus we can  discuss Borel sets in terms of their Borel codes.
To do that in an admissible set ${\mathbb{A}}$ we shall
assume that ${\mathbb{A}}$ contains some countable set
(possibly as a set of urelements).
We will say that $\omega$ is {\em realizable} in an admissible
set ${\mathbb{A}}$ if the set contains a copy of the structure
$\langle\omega,<\rangle$ as an element.
Observe that $\omega$ is realizable in any
admissible set satisfying Infinity Axiom.
On the other hand the $\omega$-model $\mathbf{HF}(\omega ,<)$
does not realize $\omega$.
Since it does not cause any misunderstanding,
we shall write $\omega$  even if we work not
with $\omega$ itself but with its copy.

\begin{definicja}\label{F1}
Let $G$ be a recursively presented group with a basis $\{V_m: m\in \omega\}$
and $X$ be a recursively presented Polish $G$-space
with a basis $\{U_n:n\in\omega\}$.
Let ${\mathbb{A}}$ be an admissible set such that
$\omega$ is realizable in ${\mathbb{A}}$.

\begin{enumerate}
\item We say that $x\in X$ is \emph{codable} in ${\mathbb{A}}$
if  the set $r_x=\{n:x\in U_n\}$ is an element of ${\mathbb{A}}$.
\item We say that the group $G$ is \emph{codable} in
${\mathbb{A}}$ if the relations $r_o=\{(k,l,n):V_kV_l\subseteq V_n\}$
and $r_i=\{(k,l):V_k^{-1}\subseteq V_l\}$ are in ${\mathbb{A}}$.
\item We say that the $G$-action on $X$ is codable in ${\mathbb{A}}$
if the relation $r_a=\{(k,i,j): V_kU_i\subseteq U_j\}$ is an element
of ${\mathbb{A}}$.
\end{enumerate}
\end{definicja}

\begin{thm}\label{nad1} 
Let ${\mathbb{A}}$ be an admissible set realizing $\omega$.
Let $G$, $X$  and $x\in X$ be codable in ${\mathbb{A}}$.
Then for any $y\in X$ if $x$ and $y$ are in the same invariant 
Borel sets with Borel codes belonging to ${\mathbb{A}}$, 
then they are in the same invariant Borel sets of Borel rank
$\alpha\le {\it{o}}({\mathbb{A}})$.
\end{thm}

{\bf Proof.} Let $r_x, r_o, r_i, r_a$ have the same meaning
as in Definition \ref{F1}.
Then easily $x$, the group operations of $G$ and the $G$-action on $X$
are recursive in $r_x, r_o, r_i$ and $r_a$ respectively.
For example, since
$\{(g,x,j): g\cdot x\in U_j, j\in\omega\}=
\bigcup\{ V_k\times U_i\times \{j\}:\ (k,i,j)\in r_a\}$,
then the $G$-action on $X$ is recursive in $r_a$.\parskip0pt

Let $\alpha< {\it{o}}({\mathbb{A}})$ be a countable ordinal.
We can find $z\in {\mathbb{A}}$ such that
$r_x, r_o, r_i$, $r_a$ are hyperarithmetic in $z$ and
$\alpha<\omega^{CK}_1(z)$.
By The Suslin-Kleene Theorem   every set which is
hyperarithmetic in $z$ has a recursive in $z$ Borel code.
Then by the assumptions
$x$ and $y$ are in the same invariant Borel sets which
are hyperarithmetic in $z$.
Hence we can apply Theorem \ref{nadef} to see that
any invariant ${\bf{\Pi}}^0_{\alpha}$-set $B\subseteq X$
containing $x$ contains $y$.
\parskip0pt

To finish the proof notice that each invariant
${\bf{\Pi}}^0_{{\it{o}}({\mathbb{A}})}$-set $B\subseteq X$
is an intersection of a family
${\it{o}}({\mathbb{A}})$-many invariant Borel sets of Borel ranks
$\xi< {\it{o}}({\mathbb{A}})$.
If $x$ belongs to such a $B$, then $y\in B$ by the first part of the proof.
$\Box$
\bigskip

At the end we shall see that Theorem \ref{nad1} indeed
generalizes the theorem of Nadel quoted in Introduction.
We deal with the translation of this theorem into
the language of the logic actions.
The set of all finitary permutations
is a dense countable subgroup which can be recursively enumerated
and turned into a recursive presentation of $S_{\infty}$.
Similarly the space $X_L$ can be recursively presented by an appropriate
recursive enumeration of the set of all ultimately equal zero sequences.
Moreover the logic action on $X_L$
is a ($\Sigma^0_1$-)recursive function.\parskip0pt

Let ${\mathbb{A}}$ be an admissible set
and ${\mathfrak{M}}\in {\mathbb{A}}$ be an $L$-structure.
Then easily $S_{\infty}$, the logic action of $S_{\infty}$ on $X_L$
and $x_{\mathfrak{M}}$ are codable in  ${\mathbb{A}}$ in the sense
of Definition \ref{F1}.\parskip0pt

To every sentence $\sigma\in L_{\omega_1\omega}$
we can assign the invariant set
$B_{\sigma}=\{x_{\mathfrak{M}}\in X_L: {\mathfrak{M}}\models \sigma\}$.
Moreover for each Borel ${\bf{\Sigma}}^0_{\alpha}$-set $B$ invariant
under the logic action there is an $L_{\omega_1\omega}$-sentence
$\sigma$ of quantifier rank $\alpha$
(the theorem of Lopez-Escobar, see \cite{sami}
and \cite{kechris}, Theorem 16.8 and its proof) such that $B=B_{\sigma}$.
Finally we see that $B$ has a Borel code in  ${\mathbb{A}}$
if and only if $B= B_{\sigma}$
for some $\sigma\in L_{\mathbb{A}}$.

\bigskip

Institute of Mathematics, University of Wroc{\l}aw, \parskip0pt

pl.Grunwaldzki 2/4, 50-384 Wroc{\l}aw, Poland \parskip0pt

E-mail: ivanov@math.uni.wroc.pl

\begin{thebibliography}{99}
\bibitem{barwise}J.Barwise, Admissible Sets and Structures,
Springer-Verlag, NY, 1975.
\bibitem{becker} H.Becker, Polish group actions: Dichotomies and
generalized elementary embeddings, J. Amer. Math. Soc. \textbf{11},
397 - 449 (1998).
\bibitem{becker2} H.Becker, Topics in invariant descriptive set
theory, Annals of Pure and Appl. Logic, \textbf{111}, 145 - 184 (2001).
\bibitem{bk} H.Becker and A.Kechris, The Descriptive Set Theory
of Polish Group Actions, Cambridge University Press, Cambridge,
1996.
\bibitem{ershov}  Yu.Ershov, \emph{Definability and Computability},
Consultants Bureau, NY, 1996.
\bibitem{hjorth} G.Hjorth, Classification and Orbit Equivalence Relations.
AMS 1991
\bibitem{kechris} A.Kechris, Classical Descriptive Set Theory,
Springer-Verlag, New York, 1995.
\bibitem{MaKe} A.S. Kechris and D.A. Martin, Infinite game and effective
descriptive set theory, Analytic sets, Academic Press  403-470, 1980 
\bibitem{mos} Y.N.Moschovakis, Descriptive Set Theory, Studies in Logic
(North Holland, Amsterdam 1980) 
\bibitem{nadel} M.Nadel, Scott sentences and admissible sets,
Ann. Math. Log. \textbf{7}(1974), 269 - 294. 
\bibitem{sami} R.Sami, Polish group actions and the Vaught
conjecture, Trans. Amer. Math. Soc., \textbf{341}, 335 - 353 (1994).
\end{thebibliography}
\end{document}